\documentclass[12pt]{amsart}
\usepackage{amsthm, amsmath}
\usepackage[psamsfonts]{amssymb} 
\usepackage{amscd}

\newsymbol\rxtimes 226F
\theoremstyle{plain}

\newtheorem{thm}[subsection]{Theorem}
\newtheorem{lem}[subsection]{Lemma}
\newtheorem{prop}[subsection]{Proposition}

\theoremstyle{remark}
\newtheorem{rem}[subsection]{Remark}

\def\pro{\noindent\textit{Proof} : }

\def\rond{\kern 2pt{\scriptstyle\circ}\kern 2pt}

\def\dim{\mathop{\rm dim}\nolimits}

\def\Z{\mathbf{Z}}

\def\P{\mathbf{P}}

\def\d{\mathbf{d}}

\def\operp{\,\buildrel {\scriptscriptstyle\perp}\over {\oplus}\,}

\def\iso{\vbox{\hbox to .8cm{\hfill{$\scriptstyle\sim$}\hfill}
\nointerlineskip\hbox to .8cm{{\hfill$\longrightarrow $\hfill}} }}
\def\la#1{\langle {#1}\rangle}
\def\sdir_#1^#2{\mathrel{\mathop{\kern0pt\oplus}\limits_{#1}^{#2}}}

\catcode`\@=11\def\eqalign#1{\null\,\vcenter{\openup\jot\m@th\ialign{
\strut\hfil$\displaystyle{##}$&$\displaystyle{{}##}$
&&\quad\strut$\displaystyle{##}$&$\displaystyle{{}##}$
\crcr#1\crcr}}\,}
\catcode`\@=12

\input amssym.def
\input amssym
\input xy
\xyoption{all}

\begin{document}
\title[The primitive cohomology  of a complete intersection]{The primitive cohomology lattice of a complete intersection}
\author[Arnaud Beauville]{Arnaud Beauville}
\address{Laboratoire J.-A. Dieudonn\'e\\
UMR 6621 du CNRS\\
Universit\'e de Nice\\
Parc Valrose\\
F-06108 Nice cedex 2, France}
\email{arnaud.beauville@unice.fr}
\date{\today}
 
\begin{abstract}
We describe the primitive cohomology lattice of a smooth even-dimensional complete intersection in projective space.
\end{abstract}

\maketitle

\section*{Introduction}
\par Let $X$ be a smooth complete intersection of degree $d$ and even dimension $n$ in projective space.
We describe in this note the lattice structure of the  primitive cohomology $H^n(X,\Z)_\mathrm{o}$ . Excluding the cubic surface and the intersection of two quadrics, we find 
\[H^n(X,\Z)_\mathrm{o} =A_{d-1}\operp k\,E_8(\pm 1)\operp \ell\, U\quad\hbox{or}\quad \langle -d\rangle \operp k'\,E_8(\pm 1)\operp \ell' \,U\]where the numbers $k,\ell,k',\ell'$ and the sign attributed to $E_8$ depend on the multidegree and dimension of $X$ -- 
see Theorem \ref{thm} for a precise statement.
 The proof is an easy consequence of  classical facts  on unimodular lattices together with the Hirzebruch formula for the Hodge numbers of $X$. 
 \par  We warn the reader  that there are many ways to write an indefinite lattice as an orthogonal sum of indecomposable ones; for instance, when $8\,|\, d$, both decompositions above hold. Still it might be useful to have a (semi-) uniform expression for this lattice. Related results, with a  different point of view, appear in \cite{LW}.
 
 \bigskip
 \section{Unimodular lattices}
\par   We will use the following standard notations for lattices: $U$ denotes the hyperbolic plane, and  $\la{d}$  the lattice $\Z\,e$ with $e^2=d$. If $L$ is a lattice, $L(-1)$ denotes the $\Z$-module $L$ with the form $x\mapsto -x^2$; if $n$ is a negative number, we put  $n\,L:=|n|\,L(-1)$.
\par  Let $L$ be an odd unimodular lattice. A primitive vector $h\in L$ is said to be \textit{characteristic} if $h\cdot x\equiv x^2$ (mod.\ 2) for all $x\in L$; this is equivalent to saying that the orthogonal lattice $h^\perp$ is even (\cite{LW}, Lemma 3.3). 

\begin{prop}\label{wall}
Let $L$ be a  unimodular lattice, of signature $(b^+,b^-)$, with $b^+,b^-\geq 2$; put $s:=b^+ -b^-$, $t=\min{(b^+,b^-)}$, $u=\min{(b^-,b^+-d)}$. Let $h$ be a primitive vector in $L$ of square $d>0$, such that $h^\perp$ is even.
\par $1)$ If $L$ is even or $8\, |\,d$	we have $h^\perp=\langle -d\rangle \operp \frac{s}{8}\,E_8\operp (t-1) \,U$. 
\par $2)$ If $L$ is odd and 	$d\leq b^+$, we have $h^\perp=A_{d-1}\operp \frac{s-d}{8}\,E_8\operp u\, U$.
\end{prop}
\pro A classical result of Wall \cite{W} tells us that $h$ is equivalent under $O(L)$ to any primitive vector $v$ of square $d$, provided $v$ is characteristic if so is $h$. If $L$ is even, we choose a hyperbolic plane $U\subset L$ with a hyperbolic basis $(e,f)$, and we put $v=e+\frac{d}{2}f$; then $v^\perp=\Z(e-\frac{d}{2}f)\operp U^\perp$, and $U^\perp$ is an indefinite unimodular lattice, hence of the form $\ p\,E_8(\pm 1)  \allowbreak \operp q\,U$. Computing  $b^+$ and $b^-$ we find the above expressions for $p$ and $q$.
\par Consider now the case when $L$ is odd. We first observe that since  $h$ is characteristic, we have $d=h^2\equiv s$ (mod. 8) (\cite{S}, V, Th. 2). Let  
$$L':=(\sdir_{i\leq d}^{\scriptscriptstyle\perp}\Z\, e_i)
\operp \frac{s-d}{8}\,E_8\operp  u\, U\quad  \hbox{with}\quad  e_1^2=\ldots =e_d^2=1\ .$$ $L'$ is odd, indefinite and has the same signature as $L$, hence is isometric to $L$. We put $v=e_1+\ldots+e_d$. The orthogonal of $v$ in $\ \sdir_{}^{\scriptscriptstyle\perp}\Z\, e_i\ $ is the root lattice $A_{d-1}$.
By Wall's theorem $h^\perp$ is isometric to $v^\perp=A_{d-1}\operp\allowbreak  \frac{s-d}{8}\,E_8\operp  u\, U$.\smallskip
\par  Suppose moreover that $8$ divides $d$, so that $8\,|\,s$. Then $L$ is isomorphic to $\Z\,e\operp\Z\,f\operp \frac{s}{8}\,E_8\operp (t-1) \,U$, with $e^2=1$, $f^2=-1$. 
Taking $v=(\frac{d}{4}+1)e+(\frac{d}{4}-1)f$ gives the result.
\qed
\bigskip
\section{Complete intersections}
\par We will check that the hypotheses of the Proposition hold for the cohomology of complete intersections; the only non trivial point is the inequality $d\leq b^+$. 
\par We will use the notations of \cite{D}. Let $\d=(d_1,\ldots ,d_c)$ be a sequence of positive integers. We denote by $V_n(\d)$ a smooth complete intersection of multidegree $\d$ in $\P^{n+c}$. We put 
\[h^{p,q}(\d)=\dim H^{p,q}(V_{p+q}(\d))\qquad\hbox{and}\qquad h^{p,q}_\mathrm{o}(\d)=h^{p,q}(\d)-\delta _{p,q}\ .\]\begin{lem}\label{leq}
$h^{p+1,q+1}(\d)\geq h^{p,q}(\d)$.
\end{lem}
\pro Following \cite{D} we introduce the formal generating series
\[H(\d)=\sum_{p,q\geq 0}h^{p,q}_\mathrm{o}(\d)y^pz^q\ \in \Z[[y,z]]\ ;\]
we define a partial order on $\Z[[y,z]]$ by writing $P\geq Q$ if $P-Q$ has non-negative coefficients. The assertion of the lemma
 is equivalent to $H(\d)\geq yzH(\d)$. The set $\mathcal{P}$ of formal series in $\Z[[y,z]]$ with this property is stable under addition and multiplication by any  $P\geq 0$ in $\Z[[y,z]]$.
The formula\[ H(d_1,\ldots ,d_c)=\sum_{{P\subset [1,d]\atop P\neq \varnothing}}\bigl[(1+y)(1+z)\bigr]^{|P|-1}\prod_{i\in P} H(d_i)\]
(\cite{D}, Cor. 2.4 (ii)) shows that it is enough to prove that $H(d)$ is in $\mathcal{P}$. \par By \cite{D}, Cor. 2.4 (i), we have
$\displaystyle H(d)=\frac{P}{1-Q}\ $ with \[P(y,z)=\sum_{i,j\geq 0}{d-1\choose i+j+1}y^iz^j\quad
\hbox{and}\quad Q(y,z)=\sum_{i,j\geq 1} {d\choose i+j}y^iz^j\ .\] 
Since $Q\geq yz$, we get  $\displaystyle\frac{1-yz}{1-Q}=1+\frac{Q-yz}{1-Q}\geq 0$, hence $(1-yz)H\geq 0$.\qed 

\begin{lem}\label{calcul} Let $d=d_1\ldots d_c$. We have:
\par  \emph{a)}  $h^{p,p}(\d)\geq d\ ;$
\par  \emph{b)}  $\ 2h^{p+1,p-1}(\d)+1\geq d$ , except in 
 the following cases:
\begin{itemize}
\item $\d=(2),(2,2)\,;$
\item $p=1$, $\d= (3),(2,3),(2,2,2),(2,2,2,2)\,;$
\item $p=2$, $\d= (2,2,2)$.
\end{itemize}
\end{lem}


\pro We first prove b) in  the case $p=1$. Then  $V_2(\d)$ is a surface $S$. The canonical bundle $K_S$  is $\mathcal{O}_S(e)$, with
$e:=d_1+\ldots d_c-c-3$; therefore $K_S^2=e^2d$. The case $e\leq 0$ is immediate, so we assume $e\geq 1$. Then  the index $K_S^2- 8\chi (\mathcal{O}_S)$ of the intersection form  is negative \cite{P};   if $e\geq 2$ we get $\chi (\mathcal{O}_S)>\frac{d}{2}$, hence $2h^{2,0}(\d)+1\geq d$.
\par  If $e=1$, we have $K_S=\mathcal{O}_S(1)$ hence $p_g=c+3$.
The possibilities for $\d$ are $(5),(2,4),(3,3), (2,2,3)$ and $(2,2,2,2)$, and we have $2(c+3)+1\geq d$ in each case except the last one. This also holds for $\d=(4)$, and the other cases are excluded.

\par  Since the index is negative, we have $h^{1,1}(\d)>2h^{2,0}(\d)+1$; this implies that a) holds (for $p=1$) except perhaps for $\d=(3),(2,2),(2,3),\allowbreak(2,2,2)$. But the corresponding $h^{1,1}$ is $7,6,19,19$, which is always $>d$.

\par    Now assume $p\geq 2$. a) follows from the previous case and Lemma \ref{leq}; similarly   it suffices to check b) for the values of $\d$ excluded in the case $p=1$. 
  Using the above formulas we find 
\[h^{3,1}(3)=1\ , \ h^{3,1}(2,3)=8\ , \ h^{3,1}(2,2,2,2)=27\ ,\ h^{4,2}(2,2,2)=6\ ,\]  
so that $2h^{p+1,p-1}(\d)+1\geq d$ for $p\geq 2$ in the three first cases and for $p\geq 3$ in the last one.\qed

\begin{thm}\label{thm}
Let $X$ be a smooth even-dimensional complete intersection in $\P^{n+c}$, of multidegree $\d=(d_1,\ldots ,d_c)$. Let $d:=d_1\ldots d_c$ be the degree of $X$, and let
$e$ be the number of integers $d_i$ which are even.
\par Let $(b^+,b^-)$ be the signature of the intersection form on $H^n(X,\Z)$; we put
$$s=b^+-b^-\quad,\quad t=\min{(b^+,b^-)}\quad,\quad u=\min{(b^+-d,b^-)}\ .$$
We assume $\d\neq (2,2)$, and $\d\neq (3),(2,2,2,2)$ when $n=2$. Then:

\begin{itemize}
\item $H^n(X,\Z)_\mathrm{o}=\langle -d\rangle \operp \frac{s}{8}\,E_8\operp (t-1) \,U\ $ if ${\frac{n}{2}+e\choose e}$ is even;
\item $H^n(X,\Z)_\mathrm{o}=A_{d-1} \operp \frac{s-d}{8}\,E_8\operp  u\,U\ $ if ${\frac{n}{2}+e\choose e}$ is odd.

\end{itemize}
\end{thm}

\par \ For a hypersurface, for instance,  we find a lattice of the form \allowbreak $\ A_{d-1} \operp p\,E_8\operp  q\,U\ $ if and only if $d$ is odd, or $d$ is even and $n$ is divisible by 4.

\medskip
\pro We  apply Proposition \ref{wall} with $L=H^n(X,\Z)$. We take for $h$ the class of a linear section of codimension $\frac{n}{2}$, so that $h^2=d$.

By \cite{LW}, Thm. 2.1 and Cor. 2.2, 
we know that
\begin{itemize}
\item $h$ is primitive;
\item $h^\perp$ is even;
\item $L$ is even or odd according to the parity of ${\frac{n}{2}+e\choose e}$. 
\end{itemize}
To apply the Proposition we only need the inequalities $b^+\geq d$ and $b^-\geq 2$. Note that the statement of the theorem holds trivially for $\d=(2)$, so we may assume $d\geq 3$.
Let us write  $n=4k+2\varepsilon $, with $\varepsilon \in\{0,1\}$.
By Hodge theory we have
$$b^+=\sum_{p+q=n\atop p\  \mathrm{even}} h^{p,q}+\varepsilon \qquad 
b^-=\sum_{p+q=n\atop p\  \mathrm{odd}} h^{p,q}-\varepsilon\ ;$$
when the inequalities a) and b) of Lemma \ref{calcul} hold this implies $b^+\geq d$ and $b^-\geq 2$, so  Proposition \ref{wall} gives the result.
\par  In the exceptional cases of Lemma \ref{calcul} b),
 the lattice $L$ is even  and we have  $b^+,b^-\geq 2$, so Proposition \ref{wall} still applies.\qed
\begin{rem}
The two first  exceptions mentioned in the theorem are well-known (\cite{D2}, Prop. 5.2): we have \allowbreak $H^2(X,\Z)_\mathrm{o}=E_6$ for a cubic surface, and $H^n(X,\Z)_\mathrm{o}=D_{n+3}$ for a $n$-dimensional intersection of two quadrics. For an intersection of 4 quadrics in $\P^6$, we have $d=16$, hence by Proposition \ref{wall} 
$$H^2(X,\Z)_\mathrm{o}=\la{-16}\operp 6\,E_8(-1)\operp 15\,U\ .$$
\end{rem}

\end{document}